\documentclass[a4paper]{article}

\usepackage[utf8]{inputenx}

\usepackage{amsmath, mathtools, amssymb, amsfonts, amsthm}
\usepackage{mathrsfs} 
\usepackage{hyperref}
\usepackage{color} 

\usepackage{style}

\title{A Jarn\'{i}k-type theorem for a problem of approximation by cubic polynomials}
\author{Alessandro Pezzoni
\footnote{Supported by EPSRC International Doctoral Scholars IDS Grant
EP/N509802/1.}}
\date{}

\begin{document}
\maketitle

\begin{abstract}
  For a given decreasing positive real function $\psi$, let $\Ap_n(\psi)$ be the
  set of real numbers for which there are infinitely many integer polynomials
  $P$ of degree up to $n$ such that $\abs{P(x)} \leq \psi(\polyheight{P})$. A
  theorem by Bernik states that $\Ap_n(\psi)$ has Hausdorff dimension
  $\frac{n+1}{w+1}$ in the special case $\psi(r) = r^{-w}$, while a theorem by
  Beresnevich, Dickinson and Velani implies that the Hausdorff measure
  $\Haus^g(\Ap_n(\psi))=\infty$ when a certain series diverges. In this paper
  we prove the convergence counterpart of this result when $P$ has bounded
  discriminant, which leads to a complete solution when $n = 3$ and $\psi(r) =
  r^{-w}$.
\end{abstract}

\section{Introduction}

Classical Diophantine Approximation studies the density of the rational numbers
in the set of real numbers, starting with Dirichlet's theorem, which
states that for any real number $x$ there are infinitely many pairs
$(p,q) \in \Z \times \Z \setminus \{0\}$ such that
\[
  \abs{q x - p} < \frac{1}{q}.
\]
This is in a sense optimal, since by Hurwitz's theorem
\[
  \abs{q \varphi - p} < \frac{1}{cq}
\]
has at most finitely many solutions $(p,q)$ as above when $c > \sqrt{5}$ and
$\varphi = (\sqrt{5}-1)/2$ (see \cite[Theorem~194]{HardyW2008} for a proof).
However, Khintchine's and Jarn\'{i}k's theorems tell us --- in a very
precise way --- how likely it is that a randomly chosen real number can be
approximated by rationals up to a certain accuracy, which is given in terms of a
decreasing function of the denominators.

More precisely, let $\psi \colon \R^+ \to \R^+$ be such a
function, called an \emph{approximation function}, and define
\[
  W_n(\psi) \defeq \left\{ x \in \R^n: \abs{q \cdot x - p} < \psi(\abs{q})
  \text{ for i.m. } (p,q) \in \Z \times \Z^n \setminus \{0\} \right\}
\]
where $\abs{q} = \max \abs{q_i}$ and ``i.m.'' is shorthand for ``infinitely
many''. Here and in what follows $\Haus^s$ denotes the usual $s$-dimensional
Hausdorff measure, which we recall in the next section.
\begin{theorem}[Jarn\'{i}k] In the above setting, for any $s \geq 0$ we have
  \[
    \Haus^s(W_1(\psi)) =
    \begin{cases}
      0 & \text{if } \sum_{q=1}^{\infty} \psi(q)^s q^{1-s} < \infty \\
      \infty & \text{otherwise}.
    \end{cases}
  \]
\end{theorem}
Furthermore, with some extra hypotheses (detailed in
corollary~\ref{cor:divergence} below) Jarn\'{i}k's theorem can be generalised to
Hausdorff $g$-measures $\Haus^g$. This gives a pretty accurate description of the
geometry of the set $W_1(\psi)$, but what if we want to consider, say,
approximations by algebraic numbers? One way of doing this is to look at the set
$\Ell_n(w)$ of real numbers $x$ such that
\[
  \abs{P(x)} < \polyheight{P}^{-w}
\]
for infinitely many integer polynomials $P$ with degree bounded above by $n$,
where $w>0$ is given and $\polyheight{P}$ denotes the \emph{height} of $P$, i.e.
the maximum absolute value of its coefficients. In other words, how close can we
get to $x$ with a root of $P$, in terms of the size of $P$?

The study of this problem dates back to Mahler. With Minkowski's linear forms
theorem one can prove that $\Ell_n(w)$ has full Lebesgue measure for any
$w \leq n$, and in 1932 Mahler conjectured that $\Ell_n(w)$ has measure $0$ for
every $w > n$ \cite{Mahle1932}. In 1969 Sprind\v{z}uk proved this in full
generality \cite{Sprin1969}, although the cases $n = 2$ and $n = 3$ had already
been settled by Kubilyus, Kasch, and Volkmann (see \cite{Volkm1961} for more
details).

The picture for Hausdorff measures, on the other hand, is a bit less clear. In
1983 Bernik proved that $\Ell_n(w)$ has Hausdorff dimension $\frac{n+1}{w+1}$
\cite{Berni1983}, and in 2006 Beresnevich, Dickinson, and Velani proved
\cite[Theorem~18]{BeresDV2006}, which specialises to the divergence part of a
Jarn\'{i}k-type theorem for Mahler's problem.  Interestingly, though, the
convergence case is not quite as straightforward as for Jarn\'{i}k's theorem and
the only results so far in this direction are for $n = 2$; these were the work
of Hussain \cite{Hussa2015} and Huang \cite{Huang2018}, who gave a more general
proof for the case of non-degenerate $\mathcal{C}^2$ plane curves.

\subsection{Hausdorff measures and dimension}

Let $X$ be a subset of $\R^n$ and let $g \colon \R^+ \to \R^+$ be a
\emph{dimension function}, i.e. a continuous, increasing function s.t.
$g(r) \to 0$ as $r \to 0$. Given $\rho > 0$, a \emph{$\rho$-cover} of $X$ is a
(possibly countable) collection $\{B_i\}$ of balls of $\R^n$ s.t. the radius
$r(B_i)$ of each ball $B_i$ lies in $(0, \rho]$ and $X \subseteq \bigcup B_i$.
Now define
\[
  \Haus^g_{\rho}(X) \defeq \inf \left\{ \sum g\left(r(B_i)\right) :
    \{B_i\} \text{ is a $\rho$-cover of } X \right\}
\]
and note that this is increasing as $\rho \to 0$. Therefore the limit
\[
  \Haus^g(X) \defeq \lim_{\rho \to 0^+} \Haus^g_{\rho}(X)
  = \sup_{\rho > 0} \Haus^g_{\rho}(X)
\]
exists, and is called the \emph{Hausdorff $g$-measure} of $X$. When
$g(r) = r^s$ for some $s \geq 0$, it is customary to write $\Haus^s(X)$ for
$\Haus^g(X)$, which is then called the \emph{$s$-dimensional Hausdorff measure}
of $X$. Moreover, if $s$ is an integer, then $\Haus^s$ is just a constant
multiple of the Lebesgue measure on $\R^s$.

If $h,g$ are two dimension functions, a straightforward standard argument shows
that if $h(r)/g(r) \to 0$ when $r \to 0$, then
\[
  \Haus^h(X) = 0 \text{ whenever } \Haus^g(X) < \infty.
\]
In particular, this implies that if $s > t \geq 0$, then $\Haus^s(X) = 0$ when
$\Haus^t(X) < \infty$. We can then define the \emph{Hausdorff dimension} of $X$
as
\[
  \dimh(X) \defeq \inf\left\{s \geq 0 : \Haus^s(X) = 0 \right\}
  = \sup\left\{ s \geq 0 : \Haus^s(X) = \infty \right\}.
\]
Finally, note that if $\dimh$ is an integer, then it coincides with the usual
``naive'' notion of dimension.

\subsection{Our setting}

From here on we will make heavy use of Vinogradov's notation $\ll$ and $\asymp$, where
$a \ll b$ means that $a \leq c b$ for some constant $c > 0$, while $a \asymp b$
means that both $a \ll b$ and $b \ll a$. We will also use subscripts to
emphasise the dependence of the implied constant on certain quantities; for
example, $a \ll_n b$ means that the implied constant $c$ depends on $n$. For
later convenience, define
\[
  \Polys_n \defeq \left\{ P \in \Z[X]: \deg(P) \leq n \right\}
\]
and
\[
  \Ap_n(\psi) \defeq \left\{ x \in \R:
    \abs{P(x)} \leq \psi(\polyheight{P}) \text{ for i.m. } P \in \Polys_n \right\}.
\]

\begin{theorem}[Beresnevich, Dickinson, Velani {\cite[Theorem~18]{BeresDV2006}}]
  Let $\mathcal{M}$ be a non-degenerate submanifold of $\R^n$ of dimension $m$.
  Let $\psi$ be an approximation function, and let $g$ be a dimension function
  such that $q^{-m} g(q)$ is decreasing and $q^{-m} g(q) \to \infty$ as
  $q \to 0$. Furthermore, suppose that $q^{1-m} g(q)$ is increasing. Then
  \[
    \Haus^g\left( W_n(\psi) \cap \mathcal{M} \right) = \infty
    \quad \text{if} \quad
    \sum_{q=1}^\infty g\left(\frac{\psi(q)}{q}\right) \psi(q)^{1-m} q^{m+n-1} = \infty.
  \]
  \label{thm:bdv_divergence}
\end{theorem}

\begin{corollary}
  Let $\psi$ be an approximation function and consider an increasing dimension
  function $g$ such that $q^{-1}g(q)$ is decreasing and $q^{-1}g(q) \to \infty$
  as $q \to 0$. Then
  \[
    \Haus^g\left( \Ap_n(\psi) \right) = \infty
    \quad \text{if} \quad
    \sum_{q=1}^\infty g\left( \frac{\psi(q)}{q} \right) q^n = \infty.
  \]
  \label{cor:divergence}
\end{corollary}


Now for any given $0 < \lambda \leq n-1$ consider
\[
  \Polys_n^{\lambda} \defeq \left\{ P \in \Polys_n:
    \abs{\polydisc{P}} \gg H(P)^{2(n-1-\lambda)} \right\},
\]
where the constant implied by the Vinogradov symbol is independent of $P$, and
\[
  \Ap_n^{\lambda}(\psi) \defeq \left\{ x \in R:
    \abs{P(x)} \leq \psi(\polyheight{P})
    \text{ for i.m. } P \in \Polys_n^{\lambda} \right\}.
\]
Note that the determinant of a polynomial $P$ of degree $n$ is a homogeneous
polynomial of degree $2n-2$ in the coefficients of $P$, and it's value is
bounded above by $c_n \polyheight{P}^{2n-2}$, for some constant $c_n$ that
depends only on $n$.

In this paper we will examine the case $n=3$ of the convergence equivalent of
corollary~\ref{cor:divergence} and provide a partial result for general $n$.
Namely, we will prove the following:

%

\begin{theorem}
  Let $\psi$ and $g$ as in Corollary~\ref{cor:divergence}. Then for any
  $0 < \lambda < 1$ we have that
  \[
    \Haus^g\left( \Ap_n^{\lambda}(\psi) \right) = 0
    \quad \text{if} \quad
    \sum_{q=1}^\infty g\left( \frac{\psi(q)}{q} \right) q^n < \infty.
  \]
  \label{thm:bounded_derivative}
\end{theorem}

For the counterpart, let's set the notation
\begin{gather*}
  \Polys_{n,\lambda} \defeq \left\{ P \in \Polys_n:
    \abs{\polydisc{P}} \ll H(P)^{2(n-1-\lambda)} \right\},\\
  \Ap_{n,\lambda}(\psi) \defeq \left\{ x \in R:
    \abs{P(x)} \leq \psi(\polyheight{P})
    \text{ for i.m. } P \in \Polys_{n,\lambda} \right\}.
\end{gather*}

\begin{theorem}
  Consider $\psi$ and $g$ as in Corollary~\ref{cor:divergence}. Let
  $\Polys_{3,\lambda}^*$ be the set of irreducible polynomials in
  $\Polys_{3,\lambda}$ and let $\Ap_{n,\lambda}^*(\psi)$ be the corresponding
  $\limsup$ set. Further assume that $0\leq \lambda < 9/20$. Then
  \[
    \Haus^g\left( \Ap_{3,\lambda}^*(\psi) \right) = 0
    \quad \text{if} \quad
    \sum_{q=1}^\infty g\left( \frac{\psi(q)}{q} \right) q^{3-2\lambda/3} < \infty.
  \]
  \label{thm:cubic_irreducible}
\end{theorem}

\begin{corollary}
  Suppose that $\psi(q) = q^{-w}$ for some $w> 0$ and that
  $0 \leq \lambda < 9/20$. As customary, write $\Ap_{3,\lambda}(w)$ for
  $\Ap_{3,\lambda}(\psi)$. Then
  \[
    \Haus^g\left( \Ap_{3,\lambda}(w) \right) = 0
    \quad \text{if} \quad
    \sum_{q=1}^\infty g\left( q^{-w-1} \right) q^{3-2\lambda/3} < \infty.
  \]
  \label{cor:cubic_power_approx}
\end{corollary}

As a special case, for $\lambda=0$ we recover Bernik's result for $n=3$, namely
that the Hausdorff dimension of $\mathscr{L}_3(w) = \Ap_{3,0}(w)$ is
$\frac{4}{w+1}$.

\section{A few lemmas on polynomials}

In this section we will collect some lemmas that we will use later in the paper.
Some we prove here, while others are taken from \cite{Sprin1969}, often restated
in a slightly simpler way that is enough for our purpose.

\begin{lemma}
  Let $P_1,\dotsc,P_k$ be integer polynomials. Then
  \[
    \polyheight{P_1\dotsm P_k} \asymp \polyheight{P_1}\dotsm \polyheight{P_k}
  \]
  where the implied constants depend only on the degrees of the polynomials.
  \begin{proof}
    Recall that the Mahler measure of a polynomial $P$ of degree $d$ is defined
    as
    \[
      \mahlerm{P} \defeq \abs{a_d} \prod_{i=1}^d \max\{1,\abs{\alpha_i}\}
    \]
    where $a_d$ and $\alpha_i$ are the leading coefficient and roots of $P$,
    respectively. Now, Mahler \cite{mahle1963} showed that $\mahlerm{P}$ satisfies
    \[
      \binom{d}{\floor{d/2}}^{-1} \polyheight{P}
      \leq \mahlerm{P} \leq
      \sqrt{d+1} \, \polyheight{P}.
    \]
    Hence the result follows by noting that the Mahler measure is
    multiplicative, which can easily be seen from its definition.
  \end{proof}
  \label{lm:product_height}
\end{lemma}

\begin{lemma}{\cite[Lemma~3]{KaschV1958}}
  Let $P$ be an integer polynomial of degree at most $n \geq 2$ and with
  non-zero discriminant. If $\alpha$ is a root of $P$, then
  \[
    \abs{P'(\alpha)} \gg \abs{\polydisc{P}}^{\frac{1}{2}} \polyheight{P}^{-n+2}
  \]
  where the implied constant depends only on $n$.
  \label{lm:derivative_and_discriminant}
\end{lemma}

\begin{lemma}{\cite[Lemma~4]{KaschV1958}}
  Let $P$ be as in lemma~\ref{lm:derivative_and_discriminant} and consider some
  $x \in \C$. If $\alpha$ is the closest root of $P$ to $x$, then
  \[
    \abs{x - \alpha} \ll
    \polyheight{P}^{n-2} \abs{\polydisc{P}}^{- \frac{1}{2}} \abs{P(x)}
  \]
  where the implied constant depends only on $n$.
  \label{lm:distance_from_root}
\end{lemma}

\begin{lemma}
  In the setting of lemma~\ref{lm:distance_from_root} write $H$ for
  $\polyheight{P}$. Furthermore, assume that
  $x \in \left[ -\frac{1}{2}, \frac{1}{2} \right)$ and that
  $\abs{P(x)} < \psi(H)$ for some approximation function $\psi$. If
  $\abs{\polydisc{P}} \gg_n H^{2(n-1-\lambda)}$ for some
  $0 < 2\lambda < 1 - \log_H \psi(H)$, then we have
  $\abs{P'(x)} \asymp_n \abs{P'(\alpha)}$ for sufficiently large $H$.
  \label{lm:comparable_derivative}
  \begin{proof}
    First, observe that by lemma~\ref{lm:derivative_and_discriminant} we have
    \[
      \abs{P'(\alpha)} \gg_n H^{n-1-\lambda} \, H^{-n+2} = H^{1-\lambda}.
    \]
    Then note that by lemma~\ref{lm:distance_from_root} we have
    \[
      \abs{x - \alpha} \ll_n H^{n-2} \, H^{\lambda+1-n} \psi(H)
      = H^{\lambda - 1} \psi(H)
      < H^{-1/2}\psi(H)^{1/2}.
    \]
    Hence we can assume $\abs{\alpha} < 1$, since this is less than $1/2$ for
    $H$ large enough. Now, by the mean value theorem we can find some $z$
    between $x$ and $\alpha$, thus with $\abs{z} < 1$, such that
    \[
      \abs{P'(x) - P'(\alpha)} = \abs{P''(z)} \abs{x-\alpha}
      \ll_n H \, H^{\lambda - 1} \psi(H)
      = H^\lambda \psi(H).
    \]
    Finally, the hypothesis on $\lambda$ implies that
    $H^\lambda \psi(H) < H^{1-\lambda}$, therefore up to choosing $H$ large
    enough we have
    \[
      \abs{P'(x) - P'(\alpha)} < \frac{1}{2}\,\abs{P'(\alpha)},
    \]
    from which it follows that $\abs{P'(x)} \asymp_n \abs{P'(\alpha)}$, as required.
  \end{proof}
\end{lemma}

\begin{lemma}
  Fix $P \in \C[X]$ and $m \in \C$. If $P(X) = a_n X^n + \dotsb + a_1 X + a_0$,
  then the coefficients of $P(X + m) = b_n X^n + \dotsb + b_n X + b_0$ are
  \[
    b_k = \sum_{j=k}^n \binom{j}{k} a_j m^{j-k}
    \quad \text{for each} \quad
    0 \leq k \leq n.
  \]
  \begin{proof}
    We proceed by induction on $n$. If $n=1$ then
    $P(X+m) = a_1 X + a_0 + m a_1$, which agrees with the above formula. Now
    assume the lemma is true for $n-1$. Since we can write
    $P(X) = a_n X^n + Q(X)$, where $Q(X) = a_{n-1} X^{n-1} + \dotsb + a_0$, we
    have that
    \[
      P(X + m) = Q(X+m) + a_n (X+m)^n
      = Q(X+m) + a_n \sum_{i=0}^n \binom{n}{i} X^i m^{n-i}.
    \]
    Thus $b_n = a_n$ and, by the induction hypothesis, for $0 \leq k \leq n-1$
    \[
      b_k = \binom{n}{k} a_n m^{n-k} + \sum_{j=k}^{n-1} \binom{j}{k} a_j m^{j-k}
      = \sum_{j=k}^n \binom{j}{k} a_j m^{j-k}. \qedhere
    \]
  \end{proof}
  \label{lm:translated_polynomial_coefficients}
\end{lemma}

\begin{note}
  While we chose to state the lemma over $\C$ for simplicity, there is nothing
  specific to it in the proof, which carries over as-is for any other
  commutative ring with unity.
\end{note}

\begin{corollary}
  In the setting of lemma~\ref{lm:translated_polynomial_coefficients} we have
  \[
    \polyheight{P(X+m)} \leq \left( 1 + \abs{m} \right)^n \polyheight{P}.
  \]
  \begin{proof}
    By lemma~\ref{lm:translated_polynomial_coefficients}, for each
    $0 \leq k \leq n$ we have
    \[
      \abs{b_k} \leq \sum_{j=k}^n \binom{j}{k} \abs{a_j}\abs{m}^{j-k}
      \leq \polyheight{P} \sum_{j=k}^n \binom{j}{k} \abs{m}^{j-k}.
    \]
    Since $\left( 1 + \abs{m} \right)^n = \sum_{s=0}^n \binom{n}{s} \abs{m}^s$,
    it is enough to prove that
    \[
      \binom{n}{s} \geq \binom{s+k}{k} = \binom{s+k}{s}
    \]
    for any $0 \leq k \leq n$ and $0 \leq s \leq n-k$. On the other hand,
    $\binom{t}{s}$ is monotonic in $t$ for $t \geq s$, which can be readily seen
    from
    \[
      \binom{t+1}{s} = \frac{t+1}{t+1-s} \, \frac{t!}{(t-s)! s!}
      \geq \binom{t}{s}.
    \]
    Therefore the observation that $n \geq j = s+k$ completes the proof.
  \end{proof}
  \label{lm:translated_polynomial_height}
\end{corollary}

\section{Proof of Theorem~\ref{thm:bounded_derivative}}

Let $I \defeq \left[-\frac{1}{2}, \frac{1}{2}\right)$. We prove the result for
$\Ap_n^{\lambda}(\psi) \cap I$, and then extend it to the whole
$\Ap_n^{\lambda}(\psi)$.
Our first goal is to estimate how much each polynomial in $\Polys_n^{\lambda}$
can contribute towards $\Ap_n^{\lambda}(\psi)$. To do so, consider some
$\varepsilon > 0$ and $Q \in \N$. For a polynomial $P \in \Polys_n^{\lambda}$
with $\polyheight{P} \leq Q$ define
\[
  \sigma_\varepsilon(P) \defeq \left\{ x \in I:
    \abs{P(x)} \leq \varepsilon, \abs{P'(x)} \geq 2 \right\}.
\]
Then let $B_n(Q,\varepsilon)$ be the union of $\sigma_\varepsilon(P)$ over all
such polynomials. We will rely on the following specialisation of
\cite[Proposition~1]{Beres1999}:

\begin{lemma}
  For any $Q > 4n^2$ and any $\varepsilon < n^{-1}2^{-n-2}Q^{-n}$ we have
  \[
    \abs{B_n(Q,\varepsilon)} \leq n2^{n+2}\varepsilon Q^n.
  \]
  \label{lm:measure}
\end{lemma}

Now, let's partition $\Polys_n^{\lambda}$ into sets
\[
  \Polys_n^{\lambda}(t) \defeq \left\{ P \in \Polys_n^{\lambda}:
    2^t \leq \polyheight{P} < 2^{t+1} \right\}
\]
and observe that
\[
  \Ap_n^{\lambda}(\psi) \cap I = \limsup \gamma_\psi(P)
  = \bigcap_{t_0 = 1}^\infty \bigcup_{t = t_0}^\infty
  \bigcup_{P \in \Polys_n^{\lambda}(t)} \gamma_\psi(P)
\]
where
$\gamma_\psi(P) \defeq \left\{x \in I: \abs{P(x)} \leq \psi(\polyheight{P})\right\}$.
Then, for $t$ large enough and for any $P \in \Polys_n^{\lambda}(t)$, letting
$\varepsilon = \psi(2^t)$ we have that
$\gamma_\psi(P) \subseteq \sigma_\varepsilon(P)$, so that the sets
$\sigma_\varepsilon(P)$ form a cover of $\Ap_n^{\lambda}(\psi) \cap I$. Indeed,
$\psi(\polyheight{P}) \leq \varepsilon$ since $\psi$ is assumed to be
decreasing. Furthermore, if $\alpha$ is the root of $P$ closest to $x$, then up
to choosing $t_0$ large enough lemma~\ref{lm:comparable_derivative} ensures that
$\abs{P'(x)}$ is comparable to $\abs{P'(\alpha)}$, hence
\[
  \abs{P'(x)} \gg_n \polyheight{P}^{1-\lambda} \geq 2^{t(1-\lambda)}
\]
so $\abs{P'(x)} > 2$, again up to choosing $t_0$ large enough.


Note that each $\sigma_\varepsilon(P)$ is a union of finitely many intervals,
the number of which is bounded above by a constant that depends only on $n$. We
can't use this directly to obtain an upper bound for the Hausdorff dimension of
$\Ap_n^{\lambda}$, though, because those intervals can be arbitrarily small, and
we also don't know how many polynomials there are in each
$\Polys_n^{\lambda}(t)$. To fix this, let's consider the sets
\[
  \tilde{\sigma}_\varepsilon(P) \defeq \bigcup_{x \in \sigma_\varepsilon(P)}
    \left\{ y \in I: \abs{y-x} < 2^{-t} \varepsilon \right\}.
\]
Clearly $\sigma_\varepsilon(P) \subseteq \tilde{\sigma}_\varepsilon(P)$.
Furthermore, by the mean value theorem, for each
$y \in \tilde{\sigma}_\varepsilon(P)$ there is a $z \in I$ which lies between $y$
and the corresponding $x \in \sigma_\varepsilon(P)$ such that
\[
  \abs{P(y) - P(x)} = \abs{P'(z)} \abs{y - x}.
\]
Since $\abs{z} < 1$ we have $\abs{P'(z)} \ll_n \polyheight{P} < 2^{t+1}$, thus
\[
  \abs{P(y)} \leq \abs{P(x)} + \abs{P'(z)} \abs{y - x}
  \ll_n \abs{P(x)} + 2\varepsilon \ll \varepsilon.
\]
Now, let $c$ be the constant implied in the above inequality, so that
$\sigma_\varepsilon(P)$ is covered by intervals in $\sigma_{c\varepsilon}(P)$ of
length at least $\ell = 2^{1-t}\varepsilon$. From this we can obtain a cover
made up of intervals of length exactly $\ell$, splitting up the larger
intervals and allowing some overlap at the edges as necessary, and by
lemma~\ref{lm:measure} the polynomials in $\Polys_n^{\lambda}(t)$ contribute at
most
\[
  \frac{\abs{B_n(2^{t+1}, c\varepsilon)}}{\ell}
  \ll_n 2^{t(n+1)} \eqdef N
\]
of these intervals. To conclude, it follows that
\begin{align*}
  \Haus^g \left( \Ap_n^{\lambda} \cap I \right)
  &\ll_n \lim_{t_0 \to \infty} \sum_{t \geq t_0} g(\ell) N \\
  &= \lim_{t_0 \to \infty} \sum_{t \geq t_0}
    g\left( \frac{\psi(2^t)}{2^{t-1}} \right) 2^{t(n+1)} \\
  &\leq \lim_{t_0 \to \infty} \sum_{t \geq t_0}
    g\left( \frac{\psi(2^t)}{2^t} \right) 2^{t(n+1)} \\
  &= 0
\end{align*}
because $g$ is assumed to be increasing, $\psi$ is decreasing, and by Cauchy's
condensation test we know that
\[
  \sum_{t \geq 0} g\left( \frac{\psi(2^t)}{2^t} \right) 2^{t(n+1)} < \infty
  \quad \text{iff} \quad
  \sum_{q \geq 1} g\left( \frac{\psi(q)}{q} \right) q^n < \infty.
\]

\subsection{Extending the argument}


Fix $m \in \Z$ and consider
$x \in \left[m - \frac{1}{2}, m + \frac{1}{2}\right)$.
Then suppose that $P \in \Polys_n^{\lambda}$ is such that
$\abs{P(x)} \leq \psi(\polyheight{P})$. Now, note that $y = x - m \in I$ and let
$Q(X) = P(X + m)$, so that $Q(y) = P(x)$. Furthermore, by
lemma~\ref{lm:translated_polynomial_height} we know that
$c \polyheight{Q} \leq \polyheight{P}$, where
$c = \left( 1 + \abs{m} \right)^{-n}$ is independent of $P$. Therefore
$Q \in \Polys_n^{\lambda}$ and
\[
  \abs{P(y)} \leq \psi\left( \polyheight{P} \right)
  \leq \psi\left( c \polyheight{Q} \right).
\]
Hence the following lemma, together with the previous argument, is enough to
complete the proof.

\begin{lemma}
  Let $0 < c_1 < c_2$. Then
  \[
    \sum_{q=1}^\infty g\left( \frac{\psi(c_1 q)}{q} \right) q^n < \infty
    \quad \text{iff} \quad
    \sum_{q=1}^\infty g\left( \frac{\psi(c_2 q)}{q} \right) q^n < \infty.
  \]
  \begin{proof}
    To begin with, assume that the series with $c_1$ converges. Since $\psi$ is
    decreasing we have $\psi(c_1 q) \geq \psi(c_2 q)$, and since $g$ is
    increasing it follows that
    \[
      \sum_{q=1}^\infty g\left( \frac{\psi(c_2 q)}{q} \right) q^n
      \leq \sum_{q=1}^\infty g\left( \frac{\psi(c_1 q)}{q} \right) q^n < \infty.
    \]
    For the other implication, note that the series with $c_1$ converges if and
    only if so does the integral
    \begin{align*}
      \int_1^\infty g\left( \frac{\psi(c_1 t)}{t} \right) t^n \, \mathrm{d}t
      &= \int_{c_1c_2^{-1}}^\infty g\left( \frac{\psi(c_2 r)}{c_2 c_1^{-1} r} \right)
         \left(c_2c_1^{-1}\right)^n r^n c_2c_1^{-1} \, \mathrm{d}r \\
      &\asymp \int_{c_1c_2^{-1}}^\infty g\left( \frac{\psi(c_2 r)}{c_2 c_1^{-1} r} \right)
         r^n \, \mathrm{d}r \\
      &\leq \int_{c_1c_2^{-1}}^\infty g\left( \frac{\psi(c_2 r)}{r} \right)
        r^n \, \mathrm{d}r.
    \end{align*}
    Finally, this last integral converges because so does the series with $c_2$,
    and because both $g$ and $\psi$ are monotonic on $(0,1)$, hence bounded on
    $(c_1c_2^{-1}, 1)$.
  \end{proof}
  \label{lm:comparable_series}
\end{lemma}

\section{Proof of Theorem~\ref{thm:cubic_irreducible}}

Just like in the proof of Theorem~\ref{thm:bounded_derivative} we will focus on
$\Ap_{n,\lambda}^*(\psi) \cap I$, after which the result immediately
extends to the whole $\Ap_{n,\lambda}^*(\psi)$. Similarly to what we did there,
let's define
\[
  \Polys_{3,\lambda}^*(t) \defeq
  \left\{ P \in \Polys_{3,\lambda}^* : 2^t \leq \polyheight{P} < 2^{t+1} \right\}.
\]
Now suppose that $P \in \Polys_{3,\lambda}^*$ and let $\sigma(P)$ be the set of
$x \in I$ such that $\abs{P(x)} \leq \psi(\polyheight{P})$. Furthermore,
let $\sigma(t)$ be the union of $\sigma(P)$ over all $P$ in
$\Polys_{3,\lambda}^*(t)$. Then, by Lemma~\ref{lm:distance_from_root}, we
know that
\[
  \abs{x-\alpha}
  \leq c \polyheight{P} \, \abs{\polydisc{P}}^{-1/2} \psi(\polyheight{P})
  \eqdef r(P,\psi)
\]
where $\alpha$ is the root of $P$ closest to $x$ and where the constant $c > 0$
is independent of $P$ and $x$. Hence $\sigma(P)$ is covered by at most three
intervals of radius $r(P,\psi)$ centred at the roots of $P$. Then
\[
  \Ap_{3,\lambda}^*(\psi) \cap I \subseteq \limsup \sigma(t)
  = \bigcap_{t_0=0}^\infty \bigcup_{t=t_0}^\infty
    \sigma(t)
\]
and
\begin{align*}
  \abs{\sigma(t)}
  &\leq \sum_{P \in \Polys_{3,\lambda}^*(t)} \abs{\sigma(P)}\\
  &\ll \sum_{P \in \Polys_{3,\lambda}^*(t)}
       \polyheight{P} \, \abs{\polydisc{P}}^{-1/2} \psi(\polyheight{P})\\
  &\ll 2^t \psi(2^t)
       \sum_{P \in \Polys_{3,\lambda}^*(t)} \abs{\polydisc{P}}^{-1/2}\\
  &\ll 2^{t(3-2\lambda/3)}\psi(2^t)
\end{align*}
because from \cite[Corollary~2]{KaliaGK2014} it follows immediately that
\[
  \sum_{P \in \Polys_{3,\lambda}^*(t)} \abs{\polydisc{P}}^{-1/2}
  \asymp 2^{t(2-2\lambda/3)}
\]
where the implied constants are absolute. Just like we did in the proof of
Theorem~\ref{thm:bounded_derivative}, let's consider a slight enlargement of
$\sigma(P)$
\[
  \tilde{\sigma}(P) = \bigcup_{x \in \sigma(P)}
                      \left\{ y \in \R : \abs{y - x} < \psi(2^t)/2^t \right\}
\]
so that for any $y \in \tilde{\sigma}(P)$ we have
\[
  \abs{P(y)} \leq \abs{P(x)} + \abs{P'(z)} \abs{x - y}
  \ll \psi(\polyheight{P}).
\]
Thus $\sigma(P) \subseteq \tilde{\sigma}(P)$ and
$\abs{\tilde{\sigma}(t)} \asymp \abs{\sigma(t)}$. It follows that we can cover
$\sigma(t)$ with at most
\[
  N \defeq \frac{\abs{\tilde{\sigma}(t)}}{\ell}
  \ll 2^{t(3-2\lambda/3)}\psi(2^t) \, \frac{2^t}{\psi(2^t)}
  = 2^{t(4-2\lambda/3)}
\]
intervals of length $\ell \defeq \psi(2^t)/2^t$. Finally, this implies that
\[
  \Haus^g\left( \Ap_{3,\lambda}(\psi) \cap I \right)
  \ll \lim_{t_0 \to \infty} \sum_{t = t_0}^\infty
      g\left( \frac{\psi(2^t)}{2^t} \right) 2^{t(4-2\lambda/3)}
  = 0
\]
since by Cauchy's condensation test we know that
\[
  \sum_{t = 0}^\infty
  g\left( \frac{\psi(2^t)}{2^t} \right) 2^{t(4-2\lambda/3)} < \infty
  \quad \text{iff} \quad
  \sum_{q = 1}^\infty
  g\left( \frac{\psi(q)}{q} \right) q^{3-2\lambda/3} < \infty.
\]

\section{Proof of Corollary~\ref{cor:cubic_power_approx}}

By Theorem~\ref{thm:cubic_irreducible} it is enough to focus on reducible
polynomials, i.e. on
$\mathcal{B} \defeq \Ap_{3,\lambda}(w) \setminus \Ap_{3,\lambda}^*(w)$. Now
consider $x \in \R$ such that $\abs{P(x)} \leq \polyheight{P}^{-w}$ for
infinitely many reducible cubic polynomials $P$ and write $P = P_1P_2$, with
$\deg(P_i) = i$. Then note that if, say,
$\abs{P_1(x)} \leq \polyheight{P_1}^{-w}$ for at most finitely many $P_1$, then
$\abs{P_1(x)} \gg \polyheight{P_1}^{-w}$ for all $P_1$ and by
lemma~\ref{lm:product_height} we have
\[
  \polyheight{P_1}^{-w} \abs{P_2(x)} \ll \abs{P_1(x) P_2(x)}
  \leq \polyheight{P}^{-w} \ll \polyheight{P_1}^{-w} \polyheight{P_2}^{-w}.
\]
It follows that for at least one $i \in \{1,2\}$ we can find a constant
$c_i > 0$ such that $\abs{P_i(x)} \leq c_i \polyheight{P_i}^{-w}$ for infinitely
many $P_i$. In other words, we have that
\[
  \mathcal{B} \subseteq \Ap_1(c_1 q^{-w}) \cup \Ap_2(c_2 q^{-w}).
\]
Similarly, by noticing that a quadratic polynomial is either irreducible or a
product or two linear polynomials, we can also find constants $c_i' > 0$ such
that
\[
  \Ap_2(c_2 q^{-2}) \subseteq \Ap_1(c_1' q^{-w}) \cup \Ap_2^*(c_2' q^{-w}),
\]
where $\Ap_2^* = \Ap_{2,0}^*$. Furthermore, without loss of generality we may
assume that $c_1 \geq c_1'$, so that
\[
  \mathcal{B} \subseteq \Ap_1(c_1 q^{-w}) \cup \Ap_2^*(c_2' q^{-w}).
\]
Then Jarn\'{i}k's theorem implies that
\[
  \Haus^g\left( \Ap_1(c_1q^{-w}) \right) = 0
  \quad \text{if} \quad
  \sum_{q=1}^\infty g(c_1 q^{-w-1})q < \infty
\]
and the proof of case II of \cite{Hussa2015} implies that
\[
  \Haus^g\left( \Ap_2^*(c_2'q^{-w}) \right) = 0
  \quad \text{if} \quad
  \sum_{q=1}^\infty g(c_2' q^{-w-1})q^2 < \infty.
\]
Finally, note that those two series converge when
\[
  \sum_{q=1}^\infty g(q^{-w-1})q^{2+\varepsilon} < \infty
\]
converges for some $\varepsilon > 0$, by the following lemma.

\begin{lemma}
  Given $c,h,\varepsilon > 0$, we have that
  \[
    S_c = \sum_{q=1}^\infty g(cq^{-w-1})q^h < \infty
    \quad \text{if} \quad
    S = \sum_{q=1}^\infty g(q^{-w-1})q^{h+\varepsilon} < \infty.
  \]
  \begin{proof}
    If $c \leq 1$ the statement follows immediately by the comparison test
    because $g$ is increasing, so let's assume $c > 1$. Then note that for any
    $\theta > 0$ we can find a $q_0 = q_0(\theta)$ such that
    \[
      c q^{-w-1} < q^{-w-1+\theta}
    \]
    for every $q > q_0$. Thus $S_c$ converges if and only if so does the
    integral
    \[
      \int_{q_0}^\infty g(cq^{-w-1})q^h \, \mathrm{d}q
      < \int_{q_0}^\infty g(q^{-w-1+\theta})q^h \, \mathrm{d}q
      = \int_{p_0}^\infty g(p^{-w-1})p^{h+\delta(h+1)} \, \mathrm{d}p
    \]
    where $\delta = \frac{w+1}{w+1-\theta} > 0$ and we applied the
    substitution $q = p^{1+\delta}$. Hence if we choose $\theta$ such that
    $\delta < \frac{\varepsilon}{h+1}$, which we can always do, then the
    convergence of $S$ implies the convergence of the last integral, and thus
    the convergence of $S_c$.
  \end{proof}
  \label{lm:approx_with_constant}
\end{lemma}

\section{Conclusions}

The main issue with proving a convergence result in the case of reducible
polynomials for more general approximation functions,
similar to what we did for Corollary~\ref{cor:cubic_power_approx}, lies in the
decoupling of the resulting inequalities
\[
  \abs{P_1(x)} \abs{P_2(x)} \leq \psi(\polyheight{P_1P_2}).
\]
Our proof carries through as-is for any other $\psi$ that is multiplicative, but
this is by no means the general case. For the case of quadratic polynomials,
Hussain \cite{Hussa2015} and Huang \cite{Huang2018} resorted to imposing a
fairly restrictive condition on the dimension function and, while this
looks artificial, in private correspondence Hussain confirmed that the
techniques used in those papers don't allow to remove it.

It would also be interesting to look into an equivalent version of
Theorem~\ref{thm:cubic_irreducible} for higher degrees, which would lead to a
complete treatment of the case of approximation functions of the form $q^{-w}$.
This will likely require some different techniques, just like Sprind\v{z}uk's
solution of Mahler's conjecture required different techniques from those
Volkmann used in his treatment of the cubic case.

\bibliographystyle{abbrv}
\bibliography{/home/bexie/documents/university/jabref/library.bib}

\bigskip

\noindent
Department of Mathematics\\
University of York\\
York, YO10 5DD\\
United Kingdom\\
E-mail: \texttt{ap1466@york.ac.uk}

\end{document}